\documentclass%% [letterpaper, 11 pt, conference,onecolumn]
              {article}

\usepackage{amsfonts,amssymb,amsmath}
\usepackage{graphicx,graphics,epsfig,xcolor}
\usepackage{cite}

 \usepackage{enumitem}

\allowdisplaybreaks
\makeatletter
\let\NAT@parse\undefined
\makeatother
\usepackage[pdftex, pdfborderstyle={/S/U/W 0}]{hyperref}

\newtheorem{remark}{Remark}
\newtheorem{assumption}{Assumption}
\newtheorem{lemma}{Lemma}

\newtheorem{theorem}{Theorem}

\def\loria{Lor{\'{\i}a}}

\title{On observer-based Asymptotic Stabilization of Non-uniformly Observable Systems via hybrid control: a Case Study}

\author{
        M. Maghenem \quad  W. Pasillas-L{\'e}pine \quad  A.  Lor{\'i}a \quad  M. Aguado-Rojas 
	\thanks{M. Maghenem is with GIPSA-Lab, University of Grenoble Alpes, CNRS, France 
          (e-mail: mohamed.maghenem@cnrs.fr). 
          W. Pasillas-L{\'e}pine is with L2S, CNRS, Univ Paris-Saclay, 
          CentraleSup\'elec, France (e-mail: wiliam.pasillas-lepine@centralesupelec.fr)
          A. Lor{\'\i}a is with L2S, CNRS, (e-mail: antonio.loria@cnrs.fr). 
          M. Aguado-Rojas is with Hitachi Astemo.
          The work of M. Maghenem and of A. Lor\'ia is supported by the French ANR via project 
          HANDY, contract number ANR-18-CE40-0010. 
}}

\begin{document}
%\bstctlcite{IEEEexample:BSTcontrol}  % this removes URL stuf, DOIs, etc from refs

\maketitle 

\begin{abstract}
For systems that are not observable at the very equilibrium of interest to be stabilized,  output-feedback stabilization is considerably challenging. In this paper we solve this control problem for the case-study of a second-order system that is bilinear and affine, both in the input and the output, but it is unobservable at the target equilibrium. The case-study is representative of a well-studied class of non-uniformly observable systems and stems from automotive control. Our main contribution is a novel certainty-equivalence hybrid controller that achieves asymptotic stabilization semiglobally. The controller relies on a switched observer that estimates the state, provided that the latter is `kept away' from the singular equilibrium. To achieve both competing tasks, stabilization and estimation, the controller also relies on the keen construction of a piecewise-constant, converging, reference. Our main results are illustrated via numerical simulations on a meaningful example.
\end{abstract}
%% \begin{keywords}
%% Non-uniformly observable systems, observers, hybrid control,   ABS.
%% \end{keywords}

\section{Introduction}

``Stabilizing a system and estimating the state are two competing processes that need to happen simultaneously in order to stabilize a partially measured system in closed loop'' \cite{brivadis_HAL22}. This statement cannot be overestimated. It is particularly meaningful when the control goal is to stabilize an equilibrium at, and near which, observability is lost. The study of such systems, for about three decades, has yielded an ever-growing literature---see, {\it e.g.}, \cite{gauthier-kupka_SIAM94,Gauthier2001deterministic,besancon_lect-notes2007,TAC-SWITCHOBS,brivadis_HAL22}. It is that beyond the significant mathematical difficulty imposed solely by designing observer design \cite{Dufour2012observer,farza-etal_AUT-2017,moreno-besancon_AUT-2021}, the stabilization problem is well-motivated by concrete engineering applications, such as sensorless motor control \cite{IbarraRojas20041079}, bioreactor systems \cite{moreno_JPC-2015,Rapaport2019robust}, electrical systems \cite{riedinger2009_7074959}, and automotive applications \cite{Hoang2014}, in which stabilization covers importance beyond intellectual curiosity. 

Despite elegant solutions for particular instances such as that studied in \cite{brivadis_CDC21-9683664}, the asymptotic output-feedback stabilization of an unobservable equilibrium for nonlinear systems remains, in general, an open question. This is the case for the representative case-study, of interest here and given by the equations\\[-12pt]
\begin{subequations}\label{abs}
     \begin{eqnarray}
         \label{abs:a} \dot{z}_1 & = & -a z_1 z_2  +  u \\
         \label{abs:b} \dot{z}_2 & = & (c z_2 + d) z_1, \quad z_1, \ z_2 \in \mathbb{R},
     \end{eqnarray}
\end{subequations}
where $a$,  $c$,  and $d>0$, $u$ is the control input, and $y= C z = z_1$ is the measured output. 
\enlargethispage{20pt}\newpage

This system belongs to the well-studied class of bilinear systems of the form $\dot x = A(y)x + Bu$---cf. \cite{besancon_lect-notes2007}, is not observable on $\{y = 0\}$ and has the interest that, up to a change of time-scale, it represents the dynamics of the so-called extended-braking stiffness  \cite{Hoang2014,AGUADOROJAS2019452}. This variable, represented by $z_2$ in \eqref{abs:b}, is a state whose regulation translates into maximizing the braking force in the antiblock braking system (ABS) of a vehicle's rolling tire \cite{Sugai1999}. Stabilizing the origin $\{(z_1,z_2)= (0,0)\}$ for \eqref{abs}, measuring only $z_1$ (which tantamounts to the linear acceleration of the tire at the wheel-ground contact point, relative to the longitudinal acceleration of the vehicle) is the open problem that we solve here. 

In line with the keen assertion of  \cite{brivadis_HAL22}, the estimation and stabilization tasks are competing processes for \eqref{abs}. On one hand, as explained in \cite{Missie-TCST-XBS}, the state estimator relies on the plant's output to provide certain persistency of excitation in order for the estimation errors to converge. This is a well-documented estimation mechanism employed for bilinear systems \cite{besancon_lect-notes2007,al:LORZAV2}. A way to excite the system is via hybrid control, by constructing a switching reference for the plant's state \cite{AGUADOROJAS2019452,Corno2012}, but this, of course, prevents stabilization. To overcome this conundrum, in this paper we propose a novel hybrid controller that relies on the exponentially convergent observer previously proposed in \cite{Hoang2014,AGUADOROJAS2019452} and on the construction of a piecewise-constant, decreasing reference (more precisely, a function taking values in a discrete set formed by the elements of a decreasing sequence). Thus, the controller is hybrid in nature and successfully stabilizes the origin semiglobally and asymptotically. Our contribution is a natural follow up on long-standing work led by the second author \cite{Hoang2014,AGUADOROJAS2019452,Missie-TCST-XBS,TAC-SWITCHOBS}, but there is no intellectual intersection with these references, devoted entirely to the observer-design problem. In addition, relatively to them, we provide here explicit exponential bounds on the observer's estimation errors.

\section{Switched-observer design}  
\label{sec:obs}
We revisit the switched observer for system \eqref{abs} originally introduced in \cite{Hoang2014} and we refine the main result therein by providing an explicit estimation of the convergence rate. This is a fundamental step for the hybrid control design,  which stands as one of our main contributions. 

Let $\hat{z} := [\hat{z}_1 \ \hat{z}_2]^\top$ denote the estimate of  $z := [z_1\ z_2]^\top$  and consider the Luenberger-type observer---see \cite{Hoang2014-TCST}, 
\begin{subequations}
\label{254}
     \begin{eqnarray}
         %\label{254a} 
         \dot{\hat{z}}_1 & = & - a z_1 \hat{z}_2 - u + k_1(z_1) z_1 (z_1 - \hat{z}_1)\\
         %\label{254b} 
         \dot{\hat{z}}_2 & = & c z_1 \hat{z}_2  + d z_1 + k_2(z_1) z_1 (z_1 - \hat{z}_1), 
     \end{eqnarray}
   \end{subequations}
where $k_1$, $k_2: \mathbb R\to \mathbb R$ are functions to be defined. Then, the dynamics of the estimation error  $ \tilde z_i := \hat z_i - z_i$ is given by 
\begin{align} \label{eqerobs}
\begin{bmatrix}  \dot{\tilde{z}}_1 \\ \dot{\tilde{z}}_2   \end{bmatrix} 
=  z_1 (t) 
\begin{bmatrix}  
- k_1(z_1(t))  & -a \\ -k_2(z_1(t))  &  c     
\end{bmatrix}
\begin{bmatrix}  
\tilde{z}_1 \\ \tilde{z}_2   
\end{bmatrix},
\end{align}
which is a linear system with state $\tilde z$, and depends on time through the measurable output trajectory $t\mapsto z_1(t)$. That is, the latter is part of a solution to \eqref{abs}, with initial conditions $(t_o,z_o)\in \mathbb R_{\geq 0}\times \mathbb R^2$, and is defined on $[t_o,t_f)$ for any $t_o \geq 0$ and some $t_f \leq \infty$. For the sake of argument, we assume in this section that $t_f = + \infty$---{\it cf.} Remark \ref{rmk:tf}.  

For the purpose of estimating $z_2$, the goal is to define $k_1$ and $k_2$ so that $\tilde z \to 0$. To that end, following \cite{Hoang2014-TCST},  we define
 $$  k_i (z_1) := 
 \left\{ 
 \begin{matrix}  
 k_i^+ & \text{if} ~ z_1 > 0     
\\
 k_i^- & \text{if} ~ z_1 < 0
 \\ 0 & \text{\ if} ~ z_1 = 0,
\end{matrix}  \right.
\quad  \quad  i \in  \{1,2\},  $$ 
so, for all $z_1\neq 0$,  the matrix on the right-hand side of \eqref{eqerobs} can only be equal to
$$
A_1 = 
\begin{bmatrix}  
- k_1^+  & - a \\ -k_2^+  &  c     
\end{bmatrix},
\quad \mbox{or} \quad
A_2 = 
\begin{bmatrix}  
 k_1^-  &  a \\ k_2^-  &  -c     
\end{bmatrix};
$$
these matrices are both Hurwitz if $k_1^+ > c$,  $k^+_2 < -\frac{c}{a} k_1^+$, $k_1^- < c$,   $k_2^- < -\frac{c}{a} k_1^-$---see \cite{AGUADOROJAS2019452}. Furthermore, the pairs $(A_1,C)$ and $(A_2,C)$ are observable and, if  $k_1^{-}  = 2c - k_1^+$, $c k_1^+ + a k_2^+ = c k_1^- + a k_2^-$, there exists a positive definite symmetric matrix $P \in \mathbb{R}^{2 \times 2}$, such that---{\it cf.} \cite{Hoang2014-TCST}     
\begin{equation}
  \label{229}  A_i^\top P + P A_i = - C^\top C \qquad \forall i \in \{1,2 \}.
\end{equation}
\indent The fact that \eqref{229} holds is significant because, {\it albeit} an appropriate change of time-scale \cite{TAC-SWITCHOBS}, the estimation error dynamics \eqref{eqerobs} may be analyzed as a switched time-invariant system of the class considered in \cite{HESTACLASALLE}.  
To better see this, let
\begin{align} \label{eqnts}
\tau := \int^t_{t_o}  | z_1(s) | ds =: f_{z_1} (t)
\end{align}
and 
$$ 
A(w_1) := \left\{  
\begin{matrix} 
A_1  &  \text{if} ~ w_1 > 0
\\
A_2  &  \text{if} ~ w_1 < 0
\\ 
0 &  \text{if} ~ w_1 =0.
\end{matrix} 
\right. 
$$
Then, Eq. \eqref{eqerobs} is equivalent to 
\begin{align} \label{eqerobs1}
\tilde{w}' := \frac{d \tilde{w}}{d \tau} =
  A(w_1(\tau)) \tilde{w}  \qquad \forall \tau \in  \text{Im}(f_{z_1}),
\end{align}
where
 Im$(\,\cdot\,)$ stands for the image of $(\,\cdot\,)$, $w_1 : \text{Im}(f_{z_1})  \rightarrow \mathbb{R}$ and  $\tilde{w} : \text{Im}(f_{z_1})  \rightarrow \mathbb{R}^2$ are given by
\begin{eqnarray}
w_1 (\tau) & := &
\left\{ \begin{matrix} 
 z_1 ( f^{-1}_{z_1} (\tau) )  & \hspace{2.3em} \text{if} ~ \text{card} (f^{-1}_{z_1} (\tau) ) = 1 \\
 0 & \hspace{2.3em} \text{otherwise},
\end{matrix}  \right.   
\\[3pt]
\tilde w (\tau) & := &
\left\{ \begin{matrix} 
 \tilde{z} ( f^{-1}_{z_1} (\tau) )  & \text{if} ~ \text{card} (f^{-1}_{z_1} (\tau) ) = 1 \\
 \tilde{z} ( \min \{f^{-1}_{z_1} (\tau)\} ) & \text{otherwise},
\end{matrix}  \right.\label{212}
\end{eqnarray}
where $\tilde z^\top := [\tilde z_1^\top\ \tilde z_2^\top]$, and $\text{card}(\cdot)$ means cardinality. That is, from Eqs. \eqref{eqnts}--\eqref{212}, we have $\tau(t_o) := 0$ and, for all initial conditions satisfying $w_1(0) = \tilde z_1(t_o)$ and  $\tilde w(0) = \tilde z(t_o)$, we have $w_1(\tau) = z_1(t)$ and $\tilde w(\tau) = \tilde z(t)$ for all $\tau \geq 0$ and $t\geq t_o$. But if $z_1(t) = 0$ the $\tau$-clock freezes while the $t$-clock goes on. That is, $\dot \tau(t) \equiv 0$ for all $t \in \mathcal T^0 := \{ t \in \mathbb R_{\geq t_o} : z_1(t) = 0\}$ and $\dot \tau(t) > 0$ for all $t \not\in \mathcal T^0$. For the $\tau$-clock, $\tau(\mathcal T^0)$ is an instant; $f^{-1}_{z_1}(\tau(\mathcal T^0))$ does not exist, so we set $w_1(\tau(\mathcal T^0)) = \tilde w(\tau(\mathcal T^0)) = 0$. Thus, since $w_1(\tau) = 0$ on a null-measure set,
$$ 
A(w_1(\tau)) \in \{A_1,A_2\}   \quad \text{for almost all} \quad \tau \in \text{Im} (f_{z_1}),
$$
and system \eqref{eqerobs1} corresponds to a linear system switching between two modes. 

  Now, because $w_1(\tau)$ and $\tilde w(\tau)$ coincide, respectively, with $z_1(t)$ and $\tilde z(t)$, the origin for \eqref{eqerobs} is asymptotically stable if and only if so is the origin for \eqref{eqerobs1}. This fact is at the basis of Lemma \ref{lem1} below. In that regard, consider the following hypothesis, which is later proved to hold by design. 
\begin{assumption}  \label{assd}  
There exist positive constants  $\tau_d$,  $\tau_s$,  $\underline{z}$, and $\bar{z}$,  and an infinite union of  disjoint  intervals $I_d$,  such that: ($i$) $ |z_1 (t)| \geq  \underline{z}$ for all $t \in I_d$, ($ii$)  $ |z_1 (t)| \leq  \bar{z}$ for all $t \in \mathbb{R}_{\geq 0} \backslash I_d$, ($iii$) the length of each connected interval in $I_d$ is no smaller than $\tau_d$,  and ($iv$) the length of each connected interval in $\mathbb{R}_{\geq 0} \backslash  I_d$ is smaller than $\tau_s$. 
\end{assumption}

\begin{lemma}  \label{lem1}
If Assumption \ref{assd} holds, there exists $\mu>0$ and, for any $P$ solving \eqref{229}, there exist $\kappa_1$ and $\kappa_2$, such that
\begin{eqnarray}
\label{380} 
|\tilde{z}(t)| \leq \kappa_1  |\tilde{z}(t_o)|  e^{- \kappa_2 \mu (t-t_o)} ~~~  \forall t \geq t_o + T,  \ t_o \geq 0.  \qquad \
\end{eqnarray}
\end{lemma}
\vskip 4pt
\begin{remark}
The proof of Lemma \ref{lem1} is constructive and is provided in Appendix \ref{applem1}.  The lemma improves over the main result in \cite{TAC-SWITCHOBS}, by establishing exponential stability and, more importantly, an explicit stability bound. This is primordial for the control design, which relies on the knowledge of the rate of decrease of the estimation errors. 
\end{remark} 

\section{Observer-based Hybrid Control Algorithm} \label{sec:algo}

Let $z^*: \mathbb R_{\geq 0}\to \mathbb R$ be a given, piecewise-constant, reference trajectory (to be defined)  for $z_1$ and consider the simple certainty-equivalence control law 
\begin{equation}
  \label{381}  u :=  a z_1 \hat{z}_2 - k z_{1e}, \quad  z_{1e} := z_1 - z^*.
\end{equation}
Then, the tracking-error dynamics corresponds to
\begin{equation}
\label{eqsys} \dot{z}_{1e}  = - (k + a  \tilde{z}_2)  z_{1e}     + a z^{*} \tilde{z}_2.
\end{equation}
This system is input-to-state stable  with respect to  $z^*$ uniformly in balls of initial conditions.  To better see this, let $R>0$ be arbitrarily fixed. Then, after \cite{Hoang2014-TCST}, let $P \in \mathbb{R}^{2 \times 2}$ be a positive definite matrix such that the time derivative of 
\begin{equation}
  \label{407}  V_{obs} (\tilde{z}) := \tilde{z}^\top P \tilde{z},  
\end{equation}
along the solutions to \eqref{eqerobs}, verifies 
\begin{equation}
  \label{397}  \dot{V}_{obs} (\tilde{z}(t)) \leq 0 \qquad \forall\, t \geq 0
\end{equation}
---see the proof of Lemma 1 in Appendix A. Therefore, for any $R>0$,  and  all  $|z(0)| \leq  R$, one can compute $\tilde{R}$ such that  $V_{obs}(\tilde{z}(0)) \leq \lambda_{\max}(P) \tilde{R}^2$. 
  In turn, we have 
\begin{align} \label{eqSobs}
 \tilde{z}(t)^2  \leq \gamma^2 \tilde{z}(0)^2 \leq   \gamma^2 \tilde{R}^2 \quad \forall t \geq 0, 
\end{align}
with $\gamma := \sqrt{\lambda_{\max}(P)/\lambda_{\min}(P)} \geq 1$. 
Then, using the function $V(z_{1e}) := (1/2) z_{1e}^2$ and setting
\begin{align} \label{eqk}
k :=   \gamma a \tilde{R}  + k',  \quad  k' > 0,
\end{align}
we see that, in view of \eqref{eqSobs}, the derivative of $V$ along the trajectories of \eqref{eqsys} satisfies
\begin{equation}
  \label{458} \dot V(z_{1e}) \leq - k'z_{1e}^2 + a |z^*|_\infty |\tilde z_2|_\infty |z_{1e}|,
\end{equation}
where $|\phi|_\infty:=\mbox{ess\,sup}_{t\geq 0}|\phi(t)|$.  
It follows from \eqref{458} that the tracking error converges provided that so do $|\tilde z_2|$  and $|z^*|$.  
On the other hand,  under Assumption \ref{assd},  for 
$|\tilde z_2|$ to converge, it is required that $|z^*(t)|$ dwells a certain amount of time separated from zero.  To achieve these antagonistic objectives,  we design a succession of cycles indexed by $i\in \{1,2,\ldots\}$,  during each of which,  $z^*$ takes values in $\left\{ -\frac{z_{in}^*}{2^i},  \frac{z_{in}^*}{2^i} \right\}$, so let
\begin{equation}\label{379} 
S^* := \bigcup^{\infty}_{i=0}   \left\{ -\frac{z_{in}^*}{2^i},  \frac{z_{in}^*}{2^i} \right\},
\end{equation}
where $z_{in}^* >0$ is fixed by design (see below).   That is, $z^*(t)$ undergoes a sequence of commutations between two constant values during each cycle (this guarantees the decrease of $|\tilde z_2|$) and the said constants decrease as the index $i$ increases. 

\noindent  \textit{\underline{Initialization step}:}  \label{initstep}
Let $z_{in}^*$ and $R>0$ be given.  Then, initially, we set $z^*(t) = z_{in}^*$ for all $t\in [0,t_1]$, where $t_1$ is to be defined, and $\hat z_2(0)$ is chosen such that \eqref{eqSobs} holds. Then, 
$$ |\tilde{z}_2(t)| \leq  \gamma  \tilde{R} 
\qquad \forall t \geq 0.  $$  
In view of  \eqref{eqSobs} and \eqref{458},  there exists  $T >0$ such that
\begin{equation*}
  |z_{1e}(t) | \leq   \frac{2a \tilde{R}}{\sqrt{k'}} |z^*(t)|   
\qquad    \forall t \geq T,   
\end{equation*}
so, by setting $k' \geq 16 a^2 \tilde{R}^2$,   it follows that 
$$  |z_{1e}(t) | \leq  z_{in}^*/2  \qquad  \forall t \geq T,   $$
$z_1(t) \in [\frac{z_{in}^*}{2}, z_{in}^*]$, and, consequently, Assumption \ref{assd} holds.  On the other hand, there exist  $\kappa_{1o}$, $\kappa_{2o}>0$ such that---see  \eqref{380}, 
\begin{align*} 
|\tilde{z}(t)| \leq  |\tilde{z}(0)| \kappa_{1o}  \exp\left( - \frac{\kappa_{2o} z_{in}^* t}{2} \right)   \qquad \forall t \geq 0,
\end{align*}
so, for any $\varepsilon>0$,  there exists $T_o \geq T > 0$ such that 
$$  
|\tilde{z}(t)| \leq g(0) (\varepsilon/\gamma)  \qquad \forall t \geq  T_o,  \quad g(0) := 1, 
$$
and, since $\gamma\geq 1$---see below \eqref{eqSobs}---we have $|\tilde z_2(t)| \leq \varepsilon$ for all $t \geq T_o$.  

\noindent \textit{\underline{First cycle}:}
From $t_1 := T_o$,   we set $z^*$ to satisfy $|z^*| = \frac{z_{in}^*}{2}$,  moreover, the tracking error $z_{1e}$ satisfies \eqref{458} with  $|\tilde{z}_2| \leq \varepsilon$. Therefore,   a limit cycle is generated by switching $z^*$ between  $-z_{in}^*/2$ and $z_{in}^*/2$ each time $\hat{z}_2(t)$ reaches $d/2c$  or $-d/2c$,  as follows: 

\begin{enumerate}[wide = 0pt, leftmargin = 1.3em,itemsep=3pt]
\item If $\hat{z}_2(t_1) \leq 0$,  $z^*(t_1)$ is set to $\frac{z_{in}^*}{2}$.   Then,  at $t'_1 \geq t_1$ such that  $\hat{z}_2(t'_1) = \frac{d}{2c}$, which means that $z_2(t'_1) \in [\frac{d}{2c} - \varepsilon, \frac{d}{2c} + \varepsilon]$,  the reference  $z^*$ is set to $z^*(t'_1)  = - \frac{z_{in}^*}{2}$.    Then,  at $t''_1 \geq t'_1$ such that $\hat{z}_2(t''_1) = -\frac{d}{2c}$,  which  means that $z_2(t''_1) \in [-\varepsilon-\frac{d}{2c}, \varepsilon-\frac{d}{2c}]$,  the reference $z^*$ is set back to 
$\frac{z_{in}^*}{2}$.    
 
\item If $\hat{z}_2(t_1) \geq 0$, the reference is set to $z^*(t_1) = - \frac{z_{in}^*}{2}$ and the same switching rules as above apply {\it mutatis mutandis}.
\end{enumerate}

Along the first cycle,   Assumption \ref{assd} holds on $[t_1,  + \infty)$; thus,  there exist positive constants $(\kappa_{11}, \kappa_{21})$ such that 
\begin{align*}
|\tilde{z}(t)|  \leq \kappa_{11} |\tilde{z}(t_1)|  e^{- \kappa_{21}  \int^t_{t_1}  | z_1(s) | ds} \quad \forall t \geq t_1. 
\end{align*}
The first cycle  ends at $t_2 := T_o + T_1 > 0$, such that  
$$ |\hat{z}_2(t_2)| \leq d/2c,  ~~   | \tilde{z}(t_2)|  \leq  (\varepsilon/\gamma) g(1),  ~ g(1) \in (0, g(0)),  $$
{\it idem} for each succeeding cycle indexed $i \geq 2$.

\noindent \textit{\underline{$i$th cycle}:} Let $\{g(i)\}^{\infty}_{i=0} \subset (0,1)$, $g(0) := 1$ be a decreasing sequence.   From $t_i = T_o + T_1 +  \cdots  + T_{i-1}$,  the reference $z^*$ is set to satisfy $|z^*| = \frac{z_{in}^*}{2^i}$, and $z_{1e}$ satisfies \eqref{458} with  $|\tilde{z}_2| \leq \varepsilon g(i-1)$ for some $g(i-1) \in (0,1)$.  Hence, a limit cycle is generated by making $z^*$ switch between $-\frac{z_{in}^*}{2^i}$ and $\frac{z_{in}^*}{2^i}$ each time $\hat{z}_2(t) = d/2^ic$ or $\hat{z}_2(t) =-d/2^ic$,  as follows:

\begin{enumerate}[wide = 0pt, leftmargin = 1.1em,itemsep=3pt]
\item  If $\hat{z}_2(t_i) \leq 0$,  $z^*(t_i)$ is set to $\frac{z_{in}^*}{2^i}$.   Then,  at $t'_i \geq t_i$
such that  
$\hat{z}_2(t'_i) = \frac{d}{2^i c}$, which means that $z_2(t'_i) \in 
[\frac{d}{2^i c} - \varepsilon, \frac{d}{2^i c} + \varepsilon]$,  the reference  $z^*$ is set to $z^*(t'_i)  = - \frac{z_{in}^*}{2^i}$.    Then,  at $t''_i \geq t'_i$ such that $\hat{z}_2(t''_i) = -\frac{d}{2^i c}$,  which  means that $z_2(t''_i) \in [-\frac{d}{2^i c} - \varepsilon , -\frac{d}{2^i c} + \varepsilon ]$,  the reference $z^*$ is set to $\frac{z_{in}^*}{2^i}$.   

\item If $\hat{z}_2(t_i) \geq 0$, the reference is  set to $z^*(t_i) = - \frac{z_{in}^*}{2^i}$, etc.
\end{enumerate}

During the $i$th cycle,   Assumption \ref{assd} holds on $[t_i,  + \infty)$,  so there exist positive constants $(\kappa_{1i}, \kappa_{2i})$ such that 
\begin{align*}
|\tilde{z}(t)|  \leq \kappa_{1i}  |\tilde{z}(t_{i})|  e^{- \kappa_{2i}   \int^t_{t_{i}}  | z_1(s) | ds}   \quad  \forall t \geq t_{i}.
\end{align*}
The cycle ends at $t_{i+1} := T_o + T_1 + ...  + T_{i} > 0$,  such that 
\begin{align*}  
 |\hat{z}_2(t_{i+1})| \leq d/2^i c ~ \text{and} ~  
 | \tilde{z}(t_{i+1})|  \leq  (\varepsilon/\gamma) g(i)  \quad \forall t \geq t_{i+1},  
\end{align*}
A new cycle starts over and so on.

\begin{remark}\label{rmk:tf}
 Inequality \eqref{458}, all of the analysis above, and consequently that in Section \ref{sec:obs}, only hold on the maximal interval of solutions---say on $[t_o,t_f)$ with $t_f \leq \infty$. To show that $t_f =+\infty$ we assume otherwise. Then, we replace  $|\tilde z_2|_\infty$ with $|\tilde z_2|_{t_f} =: c$ in \eqref{458}, so we have $\dot V(z_{1e}(t)) \leq ac |z^*|_\infty V(z_{1e}(t))$ for all $t$ such that $|z_{1e}(t)| \geq 2$. That is, as $t\to t_f$ we have $|z_{1e}(t)| \to \infty$ and $V(z_{1e}(t))\to \infty$, but integrating on both sides of $\dot V(z_{1e}(t)) \leq ac |z^*|_\infty V(z_{1e}(t))$,  we obtain $+\infty = ac |z^*|_\infty [t_f - t_o]$, which contradicts $t_f < +\infty$.
\end{remark}

\section{Main Statement} \label{sec:main}

To analyze formally the stability of the closed-loop system composed of the plant \eqref{abs}, the controller \eqref{381}, and the observer \eqref{254}, we rely on expressing it as a hybrid system that consists in the combination of a constrained differential and a constrained difference equations, as per in \cite{goebel2012hybrid}, {\it i.e.},
\begin{align} \label{eq.hsys}
\mathcal H: & \left\{ 
\begin{array}{ccll}
 \dot x & \hspace{-.5em} = &\hspace{-.5em} F(x) & x \in C \\  
    x^+ &\hspace{-.5em} = &\hspace{-.5em} G(x) & x \in D, \end{array}\right.
\end{align}
where the state variable $x \in \mathcal{X} \subset \mathbb{R}^{n}$ has a continuous evolution while in the flow set $C \subset \mathcal{X}$ and it is allowed to jump if in the jump set $D \subset \mathcal{X}$.  The  continuous- and the discrete-time evolution of $x$ are governed by the flow and the jump maps $F : C \rightarrow \mathbb{R}_{\geq 0} \times \{0\} \times \mathbb{R}^2 \times \mathbb{R}^2 \times \{0\}$ and $G : D \rightarrow \mathcal{X}$, respectively.  Furthermore, the closed-loop state is defined as
\begin{align*}
x & := (\tau, i,z, \tilde{z},z^*)  \in \mathcal{X},  
\\
\mathcal{X} & := \mathbb{R}_{\geq 0} \times \mathbb{N} \times \mathbb{R} \times \left(-\frac{d}{c}, + \infty \right) \times \mathbb{R}^2 \times S^*. 
\end{align*}
 Then,  the  jump and flow sets are defined as follows.  
The flow set $C := {\mathop{\rm cl}\nolimits}\left( \mathcal{X} \backslash D  \right)$,   where cl$(\,\cdot\,)$ denotes {\it the closure relative to $\mathcal{X}$} and the jump set $D$ $:=  D_c \cup D_{nc}$. The set $D_c$, which determines the jump conditions {\it within} the $i$th cycle,  is 
\begin{equation}\label{Dc}
D_{c}  :=   \left\{ x \in \mathcal{X} : |\hat{z}_2|  \geq  \frac{d |z^*|}{c z_{in}^*}, ~ \hat{z}_2 z^* \geq 0  \right\}
\end{equation}
and the set $D_{nc}$, which determines the jump condition {\it from} the $i$th to the $(i+1)$th cycle, is given by
\begin{align}\label{Dnc}
\nonumber 
\hspace{-1.5em} D_{nc}  :=   
\left\{ x \in \mathcal{X} : \phantom{\frac{2}{1}} \right. \hspace{-.7em} & 
 |\hat{z}_2| \leq \frac{d  |z^*|}{c z_{in}^*},  ~  \hat{z}_2 z^* \leq 0, 
  \\ & \hspace{-3em}
|\Phi_i(\tau,0)^\top P \Phi_i(\tau,0)|^{\frac{1}{2}} \leq \lambda_{\min}(P)^{\frac{1}{2}} h(i) 
\left. \phantom{\frac{2}{1}}  \hspace{-.7em} \right\},
\end{align}
where the transition matrix $\Phi_i$ is obtained by integrating ({\it e.g.}, numerically) the equation---{\it cf.} Eq. \eqref{eqerobs1}
\begin{equation}
  \label{489}  \frac{d \Phi_i}{d \tau} =
  A(w_1(\tau + \tau_i)) \Phi_i  \qquad  \tau \geq 0,
\end{equation}
where $\tau_i := \int^{t_i}_{0}  | z_1(s) | ds$, with $t_i$ being the moment when the $i$th cycle starts, 
\begin{equation}
  \label{406} h(i) := \frac{g(i)}{g(i-1)} \in (0,1),\quad  h(0) :=  \varepsilon/(\gamma \tilde{R}),
\end{equation}
and we recall that $\{g(i)\}^{\infty}_{i=0} \subset (0,1)$, with $g(0) := 1$, is a decreasing sequence. 

The definition of the jump sets $D_c$  and $D_{nc}$ follows the rationale developed in the previous section, but certain technical aspects are also considered in order to cast the analysis in the framework of \cite{goebel2012hybrid}. The respective first inequalities in $D_c$ and $D_{nc}$ correspond to the switch conditions explained in Section \ref{sec:algo}. The constraint $\hat{z}_2 z^* \leq 0$,  which requires that the signs of $\hat{z}_2$ and  $z^*$ be different, is imposed in the definition of $D_{nc}$, while the opposite is used to define $D_c$, to render the intersection of these sets empty (the apparent intersection $\{\hat z_2 = z^* =0\}$ is void since $z^*\neq 0$ by design). Defining the jump sets $D_{nc}$ and $D_{c}$ by simply imposing a strict inequality in either set would be in better concordance with the algorithm  described in the previous section, but such definition would lead to the hybrid system being not well-posed \cite{goebel2012hybrid}.
 
The third inequality in the definition of $D_{nc}$,
\begin{equation}
  \label{600} 
  |\Phi_i(\tau,0)^\top P \Phi_i(\tau,0)| \leq \lambda_{\min}(P)h(i)^2, 
\end{equation}
 is a conservative, yet verifiable, condition that essentially tests the size of the otherwise non-measurable estimation errors $\tilde z(t) \equiv \tilde w(\tau)$. To better see this, consider the function $V_{obs}$ in \eqref{407}. Its total derivative along the solutions to \eqref{eqerobs1} satisfies $\dot V_{obs}(\tilde{w}(\tau)) \leq 0$, so $V_{obs}(\tilde{w}(\tau)) \leq V_{obs}(\tilde{w}(0))$ for all $\tau \geq 0$. Hence, equivalently, 
\[
 \tilde{w}(\tau)^\top P \tilde{w}(\tau) \leq \tilde{w}(0)^\top P \tilde{w}(0).
\]
Therefore, using the fact that $\tilde{w}(\tau) = \Phi_i(\tau,0)\tilde{w}(0)$, we see that \eqref{600} implies that, for any  $\tilde{w}(0)\in \mathbb{R}^2$,
\[
 \tilde{w}(\tau)^\top P \tilde{w}(\tau) \leq \lambda_{\min}(P)h(i)^2  \tilde{w}(0)^2,
\]
that is, $ |\tilde{w}(\tau)|^2  \leq  h(i)^2 |\tilde{w}(0)|^2$.

Then, we introduce the flow map 
\begin{align} \label{Fk}
F(x)   :=  
\begin{bmatrix}  
|z_1|  
\\  
0  
\\  
\begin{bmatrix}
  - (k + a  \tilde{z}_2)  z_{1e}     + a z^{*} \tilde{z}_2
\\
( c  z_2   + d ) z_1
\end{bmatrix}
\\  
  z_1  
\begin{bmatrix}  
- k_1(z_1)  & -a \\ -k_2(z_1)  &  c     
\end{bmatrix}
\tilde{z}
 \\ 
 0   
\end{bmatrix}.    
\end{align}
Note that in the definition of $F$ the dynamics of the discrete variables $(i,z^*)$ is null, the dynamics of $\tau$ corresponds to \eqref{eqnts}, and the dynamics of $z$ and $\tilde{z}$ are simply repeated from \mbox{\eqref{eqsys} and \eqref{eqerobs}}, respectively.  

On the other hand, the jump map is given by
\begin{equation}
  \label{Gk}
G(x) := 
\begin{bmatrix}
 \left\{
  \begin{array}{cl}
      0   & \text{if} \quad x \in D_{nc}   \\ 
   \hspace{.75em} \tau \hspace{.75em}  & \text{if} \quad x \in D_{c}   
  \end{array}
 \right.\\
 \left\{
  \begin{array}{cl}
  i + 1 & \text{if} \quad x \in D_{nc}  \\ 
  i     & \text{if} \quad x \in D_{c}
  \end{array}
 \right.\\
 \begin{array}{c}
   z \\ \tilde z
 \end{array}\\
 \left\{\hspace{.4em}
  \begin{array}{cl}
     z^*/2   & \text{if} \quad x \in D_{nc}    \\ 
     - z^*  &  \text{if} \quad x \in D_{c} 
  \end{array}
 \right.
\end{bmatrix}.
\end{equation}
 The map $G$ is designed to reset the value of $\tau$ to $0$ each time  a new  cycle starts and updates the cycle index $i$.  The variables $z$ and $\tilde{z}$ are continuous variables, so they do not change their values during jumps.  According to the algorithm previously explained, the variable $z^*$ halves its size in absolute value whenever a jump to a new cycle occurs. Otherwise, while switching within a cycle, $z^*$ only alternates sign. It is important to note that since $D_c \cap D_{nc} = \emptyset$,  then the map $G$ is continuous on $D$. This is important for the system to be well-posed \cite{goebel2012hybrid}.

In addition,  the  initial state $x_o := (\tau_o, i_o,z_o, \tilde{z}_o, z^*_o) \in \mathcal{X}$ is defined as follows. By assumption, a number $R$ is known such that $|z_o| \leq R$. Then, the estimates $\hat{z}_o$ are set so that $| \tilde{z}_o | \leq \tilde{R}$ for some $\tilde R >0$ known.   
 Hence, when a reliable estimate of $|z_o|$ is available, the Initialization step described on p. \pageref{initstep} may be skipped
by defining the initial cycle index as $ i_o := \max \{0, \kappa_1( \tilde{R} ) \}$, where 
$$ \kappa_1( \tilde{R} )  := \max \left\{ i \in \mathbb{Z} : \tilde{R}  \leq  \frac{\varepsilon g(i-1)}{\gamma}  \right\}.  $$ 
Furthermore, according to \eqref{eqnts}, $\tau_o = 0$. 
Finally, the reference trajectory $z^*$  is initialized to 
\begin{align*}
\hspace{-1cm} z^*_o  :=  
  \left\{ 
  \begin{matrix}
 \left.
\begin{matrix}
 \displaystyle\frac{z_{in}^*}{2^{i_o}} & \text{if} ~ \hat{z}_{2o}  < 0    \\[5pt]
 -\displaystyle\frac{z_{in}^*}{2^{i_o}} & \text{otherwise}
\end{matrix}
\ \right\}
& \hspace{-0.1cm} \text{if} \quad i_o \geq 1 \\[3pt]
z_{in}^*  & \hspace{-0.1cm}   \text{if} \quad i_o = 0.
\end{matrix}
\right.
\end{align*} 

Our main statement establishes semi-global attractivity of the set $\mathcal{A}  :=  \left\{ x \in \mathcal{X} : z = \tilde{z}  = 0 \right\}$ for the closed-loop system.  
That is,  for any ball of initial conditions of radius $R$, there exists a control gain $k(R)$, as defined in \eqref{eqk}, such that all trajectories converge to the set $\mathcal{A}$. In particular, the domain of attraction may be enlarged by increasing the control gain. 

\begin{theorem} \label{thm1}
Consider the closed-loop hybrid system $\mathcal H = (C,F,D,G)$ defined by \eqref{eq.hsys}--\eqref{Dnc}, \eqref{406}, \eqref{Fk}, and \eqref{Gk}.
Let $R$, $\tilde{R} > 0$ be  such that $|z_o| \leq R$ and $|\tilde{z}_o| \leq \tilde{R}$, and let $(i_o,z^*_o,\tau_o)$ be defined as above.  Then, for each $R$ and $\tilde{R}$,  we can find $k>0$ such that
\begin{itemize}
\item[(i)] each solution to\footnote{Note that $(t,j)\mapsto x(t,j)$ are defined as absolutely continuous functions mapping their hybrid domain, dom~$x \subset \mathbb{R}_{\geq 0}\times \mathbb{N}$, into $\mathbb{R}^2$. See \cite{goebel2012hybrid} for details. } $\mathcal H$ satisfies 
$ \lim_{(t+j) \rightarrow +\infty} |x(t,j)|_{\mathcal{A}} =0, $ provided that $\lim_{i \rightarrow \infty} g(i) = 0$;
\item[(ii)] 
  
there exists $\kappa \in \mathcal{K}$ and $\delta^*>0$, such that, for any $\delta \in (0,\delta^*)$, if $|\tilde{z}_o| \leq \delta$,  the system's trajectories satisfy the bound $|(z,\tilde{z})|_{\infty} \leq \kappa(|z_o|+\delta)$;  
\item[(iii)] system $\mathcal{H}$ is well posed---see \cite{goebel2012hybrid},  and its solutions are uniformly non-Zeno, that is,  there exist $T>0$ and $J \in \mathbb{N}$ such that, on any time period of length $T$,  at most $J$ jumps can occur. \hfill {\small $\square$}
\end{itemize}
\end{theorem}

Statement (i) establishes attractivity of the origin,  provided that the algorithm is initialized as shown above Theorem \ref{thm1}. In that regard, note that the initialization cycle $i_o = 0$,   described in Section \ref{sec:algo},   may be avoided if one has approximate knowledge of the initial condition $z_o$,  so as to set $\hat z_o$ 
$\delta$-close for $\delta>0$ sufficiently small. In this case,  a bound on the overshoot of the trajectories is guaranteed---see Statement (ii).   Note that this bound does not require convergence of $g$ to $0$.  

Now, in general,  the assumption in (ii) is restrictive, but not for commercial ABS systems, for which the initial condition $z_{1o}$ is often approximately known.

\subsection*{Proof of Theorem \ref{thm1}}

\underline{\textit{Proof of item (i):}}  To guarantee asymptotic convergence of $z$ to zero, we first show that Assumption \ref{assd} holds on the $i$th cycle,  for all $i \in \{1,2,...\}$.   Let  $i \in \{1,2,...\}$ be arbitrarily fixed and consider the behavior of the solutions to $\mathcal{H}$ for all $t \in \mathcal I_i$, that is, during the duration of the $i$th cycle.  Over the interval $\mathcal I_i$, the solutions to $\mathcal H$ coincide with those of the hybrid system $\mathcal{H}_i := (F_i,  G_i,  C_i, D_i)$, with  state vector 
$$ x := (z, \tilde{z},z^*)  \in \mathcal{X}_i :=  \mathbb{R}^2 \times \mathbb{R}^2 \times \left\{-\frac{z_{in}^*}{2^i}, \frac{z_{in}^*}{2^i} \right\},$$ 
 flow map 
$  F_i(x)   :=  
\begin{bmatrix}  
\begin{bmatrix}
  - (k + a \tilde{z})  z_{1e}     + a z^{*} \tilde{z}_2
\\
( c  z_2   + d ) z_1
\end{bmatrix}
\\[12pt]  
  z_1  
\begin{bmatrix}  
- k_1(z_1)  & -a \\ -k_2(z_1)  &  c     
\end{bmatrix}
\tilde{z}
 \\[8pt]   
 0   
\end{bmatrix},
$
\\
 jump map 
$G_i(x) := 
\begin{bmatrix} 
z^\top
~
\tilde{z}^\top
~ - z^*  
\end{bmatrix}^\top$, 
 jump set 
$D_{i} :=  D_{i1} \cup D_{i2}$,   
where 
\begin{align*}
D_{i1} & :=   \left\{ x \in \mathcal{X}_i : \hat{z}_2  \geq  \frac{d/c}{2^i}, ~ z^* = \frac{z_{in}^*}{2^i}  \right\},
\\ 
D_{i2}  & :=   \left\{ x \in \mathcal{X}_i : \hat{z}_2  \leq - \frac{d/c}{2^i}, ~ z^* = -\frac{z_{in}^*}{2^i}  \right\},
\end{align*}
flow set $C_i := \mathop{\rm cl}\nolimits\left( \mathcal{X}_i \backslash D_i  \right)$, and 
$ x_{o}  :=  (z_{o},   \tilde{z}_{o},  z^*_{o}),  $
such that 
\begin{equation}
\label{initcyc}
\begin{aligned}
 |\tilde{z}_{o}|  \leq \frac{\varepsilon g(i-1)}{\gamma} 
~~ \text{and} ~~  z^*_{o}  = \left\{ 
\begin{matrix}
\displaystyle \frac{z_{in}^*}{2^i} & \text{if} ~ \hat{z}_{o2}  < 0   
\\[5pt]
 - \displaystyle\frac{z_{in}^*}{2^i} & \text{otherwise}. 
\end{matrix}
\right.
\end{aligned}
\end{equation}
The solutions to $\mathcal H_i$ (and consequently to $\mathcal H$ over $\mathcal I_i$), within the $i$th cycle, jump according to the conditions defining $D_{i1}\cup D_{i2}$ and satisfy the following. 

\begin{lemma} \label{lem2-3}
Consider the hybrid system $\mathcal{H}_i$ such that \eqref{initcyc} holds and let the parameter $k$ satisfy \eqref{eqk} with $2 k' \geq  a^2 \varepsilon^2$.  Then, 
\begin{description}
\item{(i)}  the set
\begin{align} \label{eqDi}
\mathcal{D}_i := \left\{x \in \mathcal{X}_i :  |z_{1}|  \leq  
\frac{z_{in}^*}{2^{i-1}} + 
\frac{z_{in}^*}{2^{i-1}} g(i-1)   \right\} 
\end{align}  
is forward invariant and  finite-time attractive.
\item{(ii)}  Let  $x_{o} \in \mathcal{D}_i$ and let the parameter $k$ satisfy \eqref{eqk} with $2 k' \geq a^2 \varepsilon^2$. Then, there exists $T_{lmin} > 0$, independent of $i$, such that the time between each pair of consecutive jumps of the solution starting from $x_o$ is larger than $T_{lmin}$. 
\end{description}
\end{lemma}

Furthermore, after Lemma \ref{lem2-3} the following also holds (see the appendix for the proofs).

\begin{lemma} \label{lem4}
Consider the hybrid system $\mathcal{H}_i (C_i, F_i,D_i,G_i)$ such that \eqref{initcyc} holds and the parameter $k$ satisfy \eqref{eqk}. Then, for $k'$ sufficiently large and independent of $i$,  there exist positive constants  $(\tau_{di}, \tau_{si}, \bar{z}_i, \underline{z}_i)$ so that Assumption \ref{assd} holds.
\end{lemma}

By Lemma \ref{lem4}, Assumption \ref{assd} holds. Therefore, from Lemma \ref{lem1} it follows that there exist  positive constants $\kappa_{1i}$, $\kappa_{2i}$, $T_i$, and $\mu_i$ such that  
\begin{align*} 
|\tilde{z}(t)|  \leq  \kappa_{1i}  |\tilde{z}(t_{i})|  e^{- \kappa_{2i} \mu_i (t-t_{oi})}  \qquad \forall t \geq t_{i} + T_i,
\end{align*}
where $t_{i}$ is the beginning the interval $\mathcal{I}_i$.  
Hence,  in view of the second condition in \eqref{Dnc}, the interval of duration of the $i$th cycle, $\mathcal I_i$, is finite.   Now,  we use Lemmata \ref{lem2-3} and \ref{lem4} to complete the proof of Item (i) of the theorem.  We show that, for each $i_o \in \{1,2,... \}$, there exits $i^* \geq i_o$ and $t_{i^*}\in \mathcal I_{i^*}$,  {\it i.e.},  during the Cycle $i^*$, such that $x(t_{i^*})\in \mathcal D_{i^*}$.   By the definition of $\mathcal D_i$,  the convergence of $z_1(t)$ follows. 

Let $i_o \in \{1,2,... \}$ and $t_{i_o} \geq 0$ be the time at which Cycle $i_o$ starts. Assume, without loss of generality,  that $z_1(t_{i_o}) > 0$, but $x(t) \notin \mathcal{D}_{i_o}$ for all $t\in \mathcal I_i$, that $z_2(t_{i_o}) \leq 0$, and that $z^*(t_{i_o}) >0$---the same reasoning that will follow applies to any other choice of initial conditions.  For the considered choice of initial conditions,  $\hat z_2(t)$ increases until one of the following two scenarios occurs: 

\begin{enumerate}[wide = 0pt, leftmargin = 1.3em,itemsep=3pt]
\item  There exist a time instant when $\hat{z}_2(t) = \frac{d/c}{2^{i_o}}$,  in which case,  sign$(z^*(t))$ becomes negative,  so the jump to Cycle $i_o + 1$ does not occur before $\hat{z}_2(t)$ becomes,  
 again,  smaller or equal than $\frac{d/c}{2^{i_o}}$.  For this to happen,  $z_2(t)$ must decrease, that is, $z_1(t)$ must become negative---see \eqref{abs:b}---and, consequently,  $x(t)$ must enter the set $\mathcal{D}_{i_o}$.  
 
\item A jump to Cycle $i_o+1$ occurs before $\hat{z}_2$ passes $\frac{d/c}{2^{i_o}}$. In this case, either the previous scenario occurs with $i_o$ replaced by $i_o + 1$ and $x(t)$ enters $\mathcal{D}_{i_o + 1}$ within Cycle $i_o + 1$, or a jump to Cycle $i_o+2$ occurs before $\hat{z}_2$ passes $\frac{d/c}{2^{i_o + 1}}$. However,  at some point, there must exist $i^* \geq i_o$ such that $x(t)$ enters $\mathcal{D}_{i^*}$ within Cycle $i^*$.
\end{enumerate}

Next, we show that $z_2(t)$ also converges, by establishing an upperbound in the latter for all $t\in \mathcal I_{i^*}$ such that $z_1(t) \in \mathcal{D}_{i^*}$ and when  $\hat{z}_{2}(t_{i^*}) = \frac{d/c}{2^{i^*-1}}$ and $z^*(t_{i^*}) = - \frac{z_{in}^*}{2^{i^*}}$.   
The latter must happen at some point while in 
Cycle ${i^*}$.
Note that after the proof of Lemma \ref{lem4} the overshoot of $z_2(t)$ occurs during the interval $[0, T_{lmin}]$,  where $T_{lmin}$ corresponds to the time it takes $z_1(t)$ to acquire the same sign as $z^*(t)$---in this case, to becomes negative.  By virtue of the comparison Lemma, it is enough to construct a bound on the solution of 
 \begin{align*} 
\begin{matrix}
 \dot{z}_2 = \big[\!\displaystyle\max_{\ z_1 \in \mathcal D_{i^*}} \! |z_1| \,\big] [\hspace{0.5pt}cz_2 + d\,],  & z_2(0) = \frac{d/c + \varepsilon}{2^{{i^*}-1}},
 \end{matrix}
 \end{align*}
over the interval $[0, T_{lmin}]$. Clearly,  we deduce an upperbound on $z_2$ that converges to zero as $i^*$ goes to infinity. 

\underline{\textit{Proof of item (ii):}}  
By definition,  the control algorithm is initiated at Cycle $i_o$ with $i_o := \max \left\{0,  \kappa_1(\delta) \right\}$.  Furthermore, when $\delta$ is sufficiently small, we conclude that $i_o := \kappa_1(\delta)$.
Therefore,  by definition of $\kappa_1$ and \eqref{eqSobs},  we conclude that
$$  |\tilde{z}(t,j)| \leq \min 
\left\{\varepsilon g(i_o -1), \gamma \delta \right\}  \qquad \forall (t,j) \in \mathop{\rm dom}\nolimits \tilde{z}.  $$
Next,   to find an upper bound for $z_1$,  we distinguish 
between two cases: 

\begin{enumerate}[wide = 0pt, leftmargin = 1.3em,itemsep=3pt]
\item If $x_{o} \in \mathcal{D}_{i_o}$,  where  $\mathcal{D}_{i_o}$ is defined in \eqref{eqDi} and is forward invariant, then we know that there exists a class $\mathcal{K}$ function $\kappa_2$ such that 
$$ z_1(t,j) \subset \kappa_2(|\delta|) [-1, 1] \qquad \forall (t,j) \in \mathop{\rm dom}\nolimits z_1.   $$
Indeed, it is easy to see that when $\delta$ goes to zero,  $i_o$ goes to infinity, and thus $\mathcal{D}_{i_o}$ reduces to $\{ 0 \}$.

\item If $\mathcal{D}_{i_o}   \subset \{x \in \mathcal{X} : z_1 \in  
[-|z_{1o}|,  |z_{1o}|\,] \}$ we use the fact that $\mathcal{D}_{i_o}$ is finite-time attractive---see Item (i) in Lemma \ref{lem2-3}.   
Furthermore,  since the flows are unique and $z_1$ is a continuous variable,  we conclude that $[-|z_{1o}|, |z_{1o}|\, ]$ must be forward invariant.  Hence,  we obtain that
\begin{align} \label{eqze2}
|z_1(t,j)| \leq |z_{1o}| \qquad \forall (t,j) \in \mathop{\rm dom}\nolimits z_1.
\end{align} 
\end{enumerate}

Finally,  to complete the proof,  we establish an upper bound on $z_2$.  Assume,  without loss of generality,  that 
$z_{2o} > 0$ and 
$ z^*_{o} = -  2^{-i_o} z_{in}^* = 
- 2^{-\kappa_1(\delta)} z_{in}^*  =: \kappa_3(\delta).  $ 
Then, consider the following two possibilities: 

\begin{enumerate}[wide = 0pt, leftmargin = 1.3em,itemsep=3pt]
\item   If  $z_{1o} \in \mathcal{D}_{i_o}$  we conclude that the overshoot of $|z_2|$ occurs only on the interval $[0,T_{1*}]$, on which  $|z_1| \leq \kappa_2(\delta)$, and before $z_1$ becomes negative.  

\item  If $\mathcal{D}_{i_o}   \subset \{x \in \mathcal{X} : z_1 \in  
[-|z_{1o}|,  |z_{1o}|\,] \}$,  we conclude that \eqref{eqze2} holds.  Hence,  the overshoot of $z_2$ occurs only on the interval $[0,T_{2*}]$,  on which $|z_1| \leq  |z_{1o}|$,  and before $z_1$ becomes negative.  
\end{enumerate}
Thus,  after the comparison Lemma,  it suffices to assess the behavior of the solutions of 
 \begin{align*} 
  \dot{z}_2 =  \max \left\{ |z_{1o}|,  \kappa_2(|\delta|) \right\}  (cz_2 + d) ~~ \text{with} ~~ z_2(0) = |z_{2o}|, 
 \end{align*}
 over the interval $ [0, T_*]$, where $ T_* :=  \max\{T_{1*}, T_{2*}\}$ is an upper bound on the time that $|z_1|$ takes to flow from $\max \left\{ |z_{1o}|,  \kappa_2(|\delta|) \right\} $ to zero. 
 
 To complete the proof,  we show that $T_*$ can be chosen as a class $\mathcal{K}$ function of $|(z_o,\tilde{z}_o)|$.  To do so,  we use the Lyapunov function 
$ v (z_{1e}) := z_{1e}^2,  $
whose time derivative along the solutions to $\mathcal{H}_i$ satisfies 
\begin{align*}
  \dot{v} & =  - 2 k'  v + 2 a z^* \tilde{z}_2 z_{1e} 
 \leq - 2 k'  v + 2 a \varepsilon  \frac{z_{in}^*}{2^{i_o}} g(i_o-1)  |z_{1e}|.
\end{align*}
Let $k' \geq 1$. After the triangle inequality and $g(i_o-1)<1$, 
\begin{align*}
\dot{v} 
& \leq -  k'  v + \frac{z_{in}^{*2}}{2^{2i_o}}   \frac{a^2 \varepsilon^2}{4^{\kappa_1(
\delta)}}   \leq -  k'  v +   a^2  \kappa_4(\delta).
\end{align*} 
 Then, $T_*$ corresponds to the time elapsed for the solution of $\dot{v} = -  k'  v +   a^2  \kappa_4(\delta)$ to attain the value $v(t) = \frac{z_{in}^{*2}}{2^{2i_o}}$ from  $v(0) := \left[\max \left\{ |z_{1o}|, \kappa_2(\delta)  \right\} +  \frac{z_{in}^{*}}{2^{i_o}} \right]^2$. We obtain that, since $i_o:=\max\{0, \kappa_1(\delta)\}$, $T_*$ is upper bounded by a class $\mathcal{K}$ function  of $(|z_o| + \delta)$.  

\underline{\textit{Proof of item (iii):}} 
After \mbox{\cite{goebel2012hybrid}},  system $\mathcal H$ is well-posed if the sets $C$ and $D$ are closed relative to $\mathcal{X}$ and $F$ and $G$ are continuous on $C$ and $D$,  respectively.   It is easy to conclude that our closed-loop hybrid system 
$\mathcal{H}$ satisfies the hybrid basic conditions which require the sets $C$ and $D$ to be closed and the maps $F$ and $G$ to be continuous. Note that both $C$ and $D$ are closed subsets relative to $\mathcal{X}$, $F$ is smooth and $G$ is  continuous on $D = D_c \cup D_{nc}$ since both $D_c$ and $D_{nc}$ are closed relative to $\mathcal{X}$ and their intersection is empty. 

Next,  we show that the closed-loop solutions are uniformly  non-Zeno.   To do so,  we note that within a same Cycle $i$, and between each two consecutive jumps,  $\hat{z}_2(t)$  flows from $-\frac{d}{2^{i-1} c}$ to  $\frac{d}{2^{i-1} c}$ back and forth.  The latter flow phase takes a time we denote by $T_{li}$.  After Item (ii) in Lemma \ref{lem2-3}, there exists a uniform lower bound  $T_{lmin} > 0$ such that $T_{li} \geq T_{lmin}$ for all $i \in \{1,2,...\}$, provided that $z_{1o} \in \mathcal{D}_i$. In general, after Item (i) of Lemma \ref{lem2-3},  $z_1$ must reach $\mathcal{D}_i$ in finite time while in Cycle $i$; otherwise, only one jump occurs within Cycle $i$. 

On the other hand, for a jump from Cycle $i$ to Cycle $i+1$ to occur,   the variable $\tau$ must  flow so that $|\Phi_i(\tau,0)^\top P \Phi_i(\tau,0)|^{\frac{1}{2}}$ decreases from $|P|^{\frac{1}{2}}$ to $\frac{\lambda_{min}(P)^{\frac{1}{2}}}{2}$,  where $\Phi_i$ is  defined by \eqref{489}.  We show the existence of a strictly positive lower bound on the time the latter decrease process takes.  To that end, we use $V = |\tilde{w}|^2$ and the fact that 
\begin{align*}
  \tilde{w}^\top & \left[A(w_1(\tau)) + A(w_1(\tau))^\top \right]  \tilde{w}
 \geq  - \eta |\tilde{w}|^2    \quad \forall\,\tau \geq 0, \,  \tilde w \in \mathbb{R}^2,  
\end{align*}
where $\eta:= \max_{i\in \{1,2\}} \{ |A_i + A_i^\top| \}$.  Then, $ V'(\tau) \geq - a V(\tau) $ for all  $\tau \geq 0$, which, by using the comparison Lemma, implies that $V(\tau) \geq e^{-a \tau} V(0)$ and, in turn,  for each $\tilde{w}(0) \in \mathbb{R}^2$, 
$  \tilde{w}(0)^\top \Phi_i(\tau,0)^\top  \Phi_i(\tau,0) \tilde{w}(0)  \geq e^{-a \tau} |\tilde{w}(0)|^2,  $ 
hence,  since for our case $|P| > \lambda_{min}(P)$ there exists $\tau^* > 0$ such that,  for each $i \in \{1,2,... \}$, and for all $\tau \in [0,\tau^*]$,
$$ |\Phi_i(\tau,0)^\top P \Phi_i(\tau,0)|  \geq \lambda_{min}(P) \geq \lambda_{min}(P)  h(i).  $$

\section{A numerical example}

To illustrate our theoretical findings, we performed some numerical simulations on a model of extended braking stiffness dynamics, described in detail in \cite{AGUADOROJAS2019452}. 
The switching logic is determined by the conditions in \eqref{Dc} and \eqref{Dnc}, with $a=375$, $c=24$, and $d=12.5$. The scenario defined by these parameters corresponding to that of a hard braking on a dry road. We took, moreover, a gain $k=500$ for the controller and gains $k_1^+=40$ and $k_2^+=-3$ for the observer. The initial wheel acceleration reference $z^\ast_o$ was set to~$75$ and the sequence $h$ was chosen as
$$ h(i) := 
\left\{ 
\begin{array}{cl}
1/(1+ \frac{1}{4^i})  & \text{if} \quad i \in \{1,2,...,8\},
\\ 
1/2 & \text{if} \quad i \in \{9,10,...\}.
\end{array}\right. $$
This choice of $h$ results from extensive numerical tests, which showed that the natural choice $h(i) \equiv 1/2$, which satisfies the conditions of Theorem \ref{thm1}, leads to relatively slow convergence of $z$, while $h(i) = 1/(1+ \frac{1}{4^i})$ for all $i \in \mathbb{N}$ violates the condition t $\lim_{i \rightarrow \infty} g(i) = 0$, with $g(i) = \prod_1^i h(j)$, in Item (i) of Theorem \ref{thm1}. 

\begin{figure}[h!]
%  \hspace{0mm}
\centerline{\includegraphics[scale=0.47]{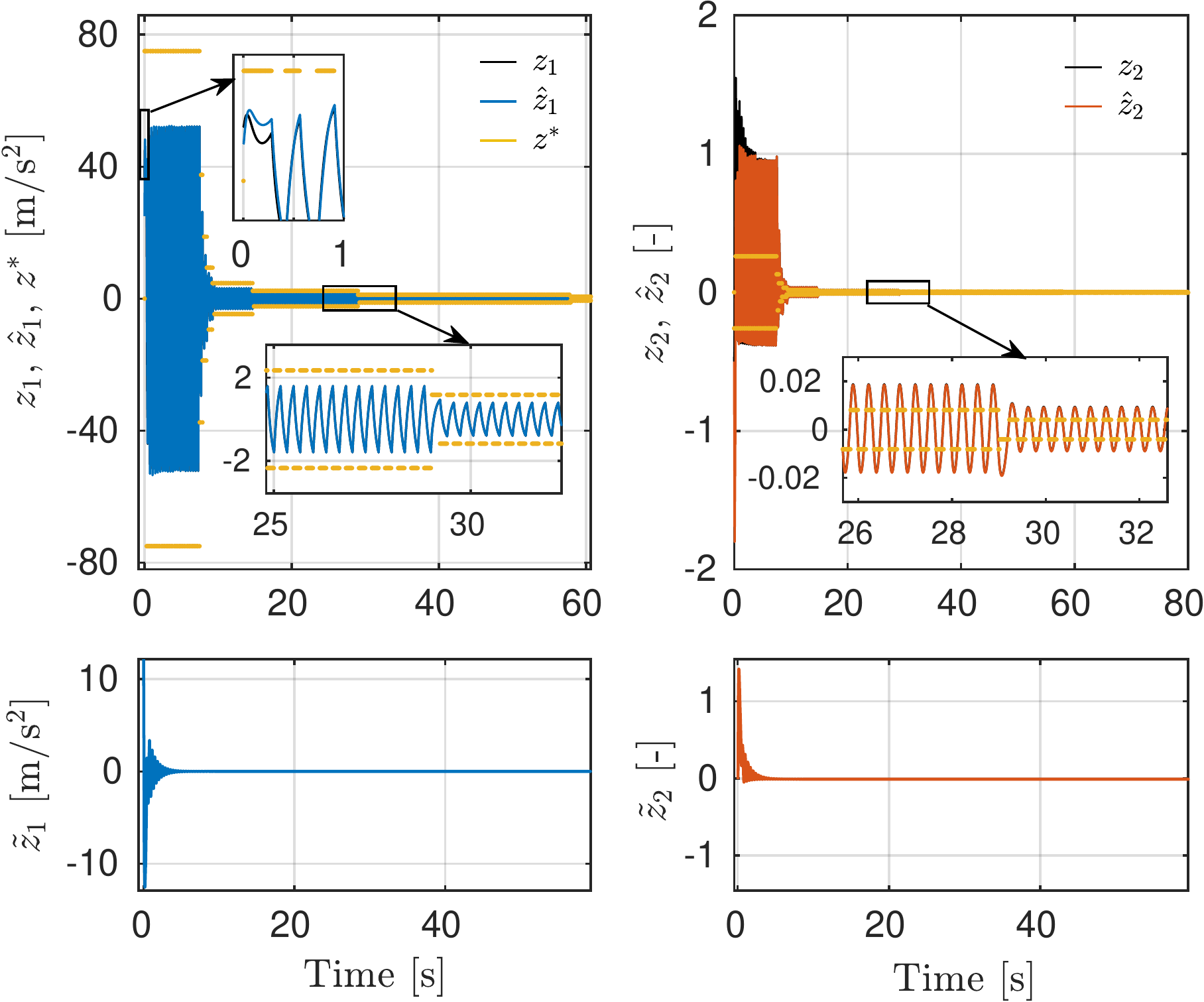}}
  \caption{System's and observer's response}
  \label{fig1}
\end{figure}
%\vskip -9pt
\begin{figure}[h!]
%\hspace{0mm}
\centerline{\includegraphics[scale=0.45]{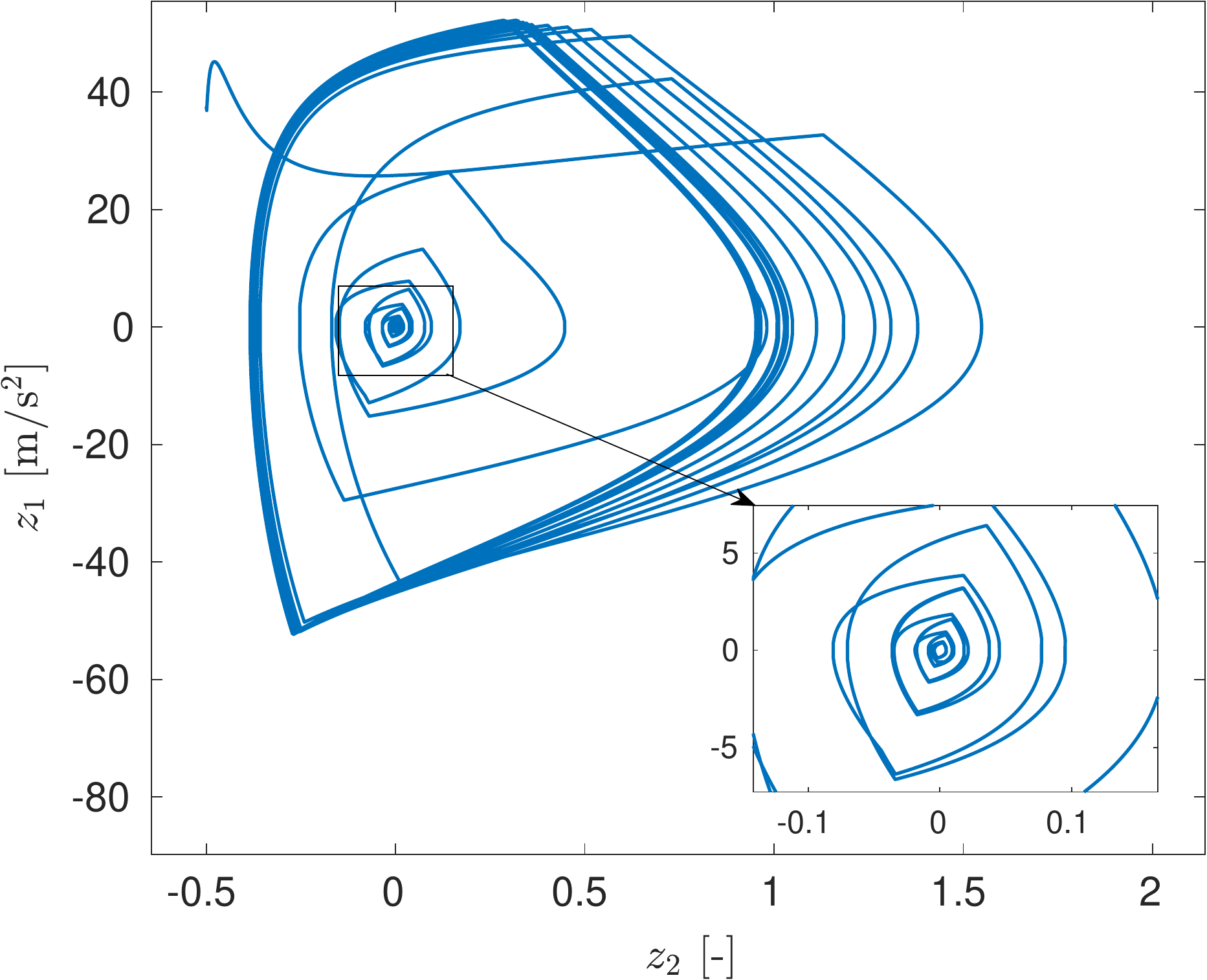}}
  \caption{System's response on the phase-plane }
    \label{fig2}
\end{figure}

The numerical results are illustrated in Figs. \ref{fig1} and \ref{fig2}. The NE plot in Fig. \ref{fig1} shows (in black) the trajectories of the measured output, $z_1$, and the corresponding estimate generated by the observer (in blue). The yellow (apparently) solid lines correspond to the piece-wise constant reference $z^*$ taking values in the discrete set $S^*$ defined in \eqref{379}, hence the ``staircase-type'' graph. It is noticeable that, {\it e.g.}, approximately $\forall \, t\in (0,8)$, $|z^*(t)| \approx 75$, while  $\forall\, t\in (15,28)$ $|z^*(t)| \approx 2.3$. Other values are discernible over other intervals. These correspond to the duration of the different cycles triggered as explained in Section \ref{sec:algo}. Within each cycle, the reference keeps switching between the positive and negative value of the same constant. In the zoomed windows, one can appreciate that, even if the system's output does not attain the imposed reference (because the latter switches) the observer keeps tracking the plant's states; the trajectories of $z_2$ and $\hat z_2$ are depicted in the NW plot. It is the ``persistent excitation'' induced by the switchings that the observer keeps up with a close estimate of the states. Then, as described in the algorithm, successive cycles, at the beginning of which the reference $z^*$ is halved, are generated. Thus, the reference, the plant's states and the observer errors, all tend to zero asymptotically.  In Fig. \ref{fig2} one may appreciate the stabilization mechanism by looking at the response on the phase plane.

\section{Discussion and Future Work}

Simultaneous estimation and stabilization at an equilibrium where observability is lost is a challenging, but not impossible task, at least for certain bilinear systems. Although in this work we focuss on  a particular plant, our results set the basis for future work oriented towards broadening the applicability of our switching-observer-based hybrid controller to other bilinear and, more generally, non-uniformly observable systems. Furthermore, an intriguing problem remains open even for the case-study aaddressed here, how to construct a smooth output-feedback observer-based controller. 

\appendix
\section{Appendix}

\subsection{Proof of Lemma \ref{lem1}}
\label{applem1}

Consider system \eqref{eqerobs1}, which is equivalent to \eqref{eqerobs}. It is a linear system that switches  between two modes defined by the matrices $A_{1}$ and $A_2$, which are both Hurwitz and the pairs $(A_1,C)$ and $(A_2,C)$ are observable. Let \cite{Hoang2014-TCST} generate a positive definite matrix  $P$ such that \eqref{229} holds.  Then, the derivative of
$$ V_{obs} (\tilde{w}) := \tilde{w}^\top P \tilde{w}  $$ 
along the solutions to the switched linear system in \eqref{eqerobs1} verifies 
$$ V'_{obs}(\tilde{w}) = - \tilde{w}^\top C^\top C  \tilde{w} \leq 0,  $$
which implies uniform global stability ({\it i.e.}, uniform stability and uniform global boundedness)  of the origin for \eqref{eqerobs1}. Furthermore, under Assumption \ref{assd},   $\text{Im} (f_{z_1}) = \mathbb{R}_{\geq 0}$ and there exist an infinite union of disjoint intervals, denoted $\bar{I}_d := f_{z_1} (I_d)$,  such that: ($i$) $|w_1 (\tau)| > \underline{z}$   for all $\tau \in \bar{I}_d$, ($ii$) the length of each connected interval in $\bar{I}_d$ is no smaller than $\tau_d \underline{z}$, and ($iii$) the length of each connected interval in $\mathbb{R}_{\geq 0} \backslash  \bar{I}_d$ is smaller than $\tau_s \bar{z}$. 

Now, inspired by \cite{HESTACLASALLE}, let $\lambda > 0$ and  $\bar{c} := e^{\max \{ |A_1|,  |A_2| \} \bar{z} \tau_s}$, and let \cite[Lemma 9]{HESTACLASALLE} generate $K_1 \in \mathbb{R}^{2}$ and $K_2 \in \mathbb{R}^{2}$ such that,  for each $i \in \{1,2\}$,  
$$ 
\big| e^{  \left(A_i + K_i C \right) \tau} \big|  \leq \frac{1}{\bar{c}}  e^{-2  \lambda (\tau - \frac{\tau_d \underline{z}}{2} )}   \quad \forall t \geq \frac{\tau_d \underline{z}}{2}.
$$

Furthermore, let $ k :=  \max \{ |K_1|, |K_2| \}$,    $\bar{k} :=  p_M \bar{c}^2 k^2 / \lambda$, $\gamma := p_M/p_m$, where $p_m I \leq P \leq p_M I$, and $\kappa_1 := \gamma [ \bar{k} + 2 \bar{c}^2 \gamma ]/[ \rho (1 + \bar{k})]$. Finally, let $L > 0$ be such that 
$$  \rho := \frac{ 2 \gamma \bar{c}^2 e^{-2 \lambda L}   + \bar{k}}{1 + \bar{k}}  < 1.  $$

According to the proof of \cite[Lemma 5]{HESTACLASALLE},  there exists a map $w_1 \mapsto K(w_1) \in \{ K_1,K_2,0 \}$ such that, along each map $\tau \mapsto w_1(\tau)$ enjoying the properties ($i$)--($iii$) listed above, %we have
$$ |\Phi_{\bar{z}_1} (\tau_1,\tau_o)|  \leq \bar{c} e^{-\lambda (\tau_1 - \tau_o)} \qquad \forall \tau_1 \geq \tau_o \geq 0,  $$
where $\Phi_{\bar{z}_1} $ is the transition matrix of the system 
$$ \tilde{w}' = \left[ A(w_1(\tau))  + K(w_1 (\tau)) C   \right] \tilde{w}.    $$
Now, after Assumption \ref{assd}, there exists $\mu>0$ such that 
$$  \int^{t}_{t_o}  | z_1(s) | d s \geq \mu(t-t_o)
  \quad \forall t \geq t_o + T,~\forall t_o\geq 0, 
$$
so, using the proof of Theorem 4 in \cite{HESTACLASALLE},  we conclude that  
$$ |\tilde{w}(\tau_1) | \leq  \kappa_1     
\rho^{\frac{\tau_1 - \tau_o }{L}}  |\tilde{w}(\tau_o)|  \qquad \forall \tau_1 \geq \tau_o \geq 0.  $$
Therefore, the statement follows defining $\kappa_2  := -\frac{\ln(\rho)}{L}$ and using the fact that $\tilde z(t) \equiv \tilde w(\tau)$.

\subsection{Proof of Lemma \ref{lem2-3}}

\subsubsection{Proof of Item (i)} 
We first use the fact that $|\tilde{z}_{o}| \leq \frac{\varepsilon g(i-1)}{\gamma}$, together with \eqref{eqSobs}, to conclude that
\begin{align*} 
|\tilde{z}(t,j)| \leq \varepsilon g(i-1)  \qquad \forall (t,j) \geq \mathop{\rm dom}\nolimits \tilde{z}.  
\end{align*}

Next, we use the Lyapunov function 
$v(z_1) := z_1^2$,
whose time derivative along the solutions to $\mathcal{H}_i$ satisfies 
\begin{align*}
\dot{v}
\leq  - k' v + \frac{2 k' {z_{in}^*}^2}{4^{i}}  + \frac{2 a^2 \varepsilon^2 {z_{in}^*}^2}{4^{i}} g(i-1)^2.
\end{align*}
As a result, we obtain 
\begin{align*}
v(t,j) & \leq  v(0,0) e^{-k' t}  + \frac{2 {z_{in}^*}^2}{2^{2i}}  \left[1 -  e^{-k' t} \right]
\\ &
\qquad + \frac{a^2 \varepsilon^2 {z_{in}^*}^2}{ (2k') 4^{ (i-1)}} g(i-1)^2 \left[1  -  e^{-k' t} \right],
\end{align*}
 so,   by choosing $k'$ such that $ k' \geq  \frac{a^2 \varepsilon^2}{2}$, we get 
\begin{align*}
v(t,j) &  \leq  v(0,0) e^{-k' t} 
\\ &
+ \left[ \frac{{z_{in}^*}^2}{2^{2(i-1)}} + \frac{{z_{in}^*}^2}{4^{(i-1)}} g(i-1)^2 \right]  \left[1 -  e^{-k' t} \right]
\\ & \leq \max \left\{v(0,0),   \left[ \frac{{z_{in}^*}^2}{2^{2(i-1)}} + 
\frac{{z_{in}^*}^2}{2^{2(i-1)}} g(i-1)^2 \right] \right\}.
\end{align*}
Hence,  the set $\mathcal{D}_i$ is finite-time attractive and forward invariant. 

\subsubsection{Proof of Item (ii)}  
 If $x_{oi} \in \mathcal{D}_i$, using the comparison lemma,  a lower bound on the time  between each two consecutive jumps of $\mathcal{H}_i$ can be obtained by computing the time the solution to 
\begin{equation}
  \label{975}  \dot{z}_2 = - \left( \frac{2 z_{in}^*}{2^{i-1}} \right)  (cz_2 + d), 
\end{equation}
with initial conditions $z_{o2} = \frac{d/c - \varepsilon}{2^{i-1}}$. That is, the time that takes to reach $\frac{-d/c + \varepsilon}{2^{i-1}}$. To compute such time,  we introduce the new time scale 
$ \tau := \left( \frac{2 z_{in}^*}{2^{i-1}}  \right) t $, to obtain in the new time scale 
$ z'_2 = -  cz_2 - d. $
By solving the latter equation, we obtain    
$$ z_2(\tau) = \sigma_i e^{-c \tau}  - \frac{d}{c}  \left[ 1 -  e^{-c \tau}   \right],
\quad \sigma_i:= \left[ \frac{d/c - \varepsilon}{2^{i-1}} \right] $$
and we use the latter to  solve $z_2(\tau) = \sigma_i$ for  $\tau$.  Reordering terms, we obtain $  e^{-c \tau} [ d/c + \sigma_i ]   = d/c  -  \sigma_i$.  Hence, 
$$ \tau  =  \frac{1}{c}    \ln   \left[ \frac{d/c + \sigma_i }{ d/c - \sigma_i }  \right]  =  
 \frac{1}{c}    \ln   \left[ 1 + 2 \frac{ \sigma_i }{ d/c -\sigma_i }   \right].   $$ 
This implies that, in the original time scale,  the length of the interval $[t_i,t_{i+1}]$ between two jumps of the solution to \eqref{975}, denoted $T_{li}$, satisfies 
 $ T_{li}   \geq    \frac{2^{i-1}}{2 c z_{in}^*}    \ln   \left[ 1 + 2 \frac{ \sigma_i }{  d/c - \sigma_i }   \right],   $
which is separated from zero, {\it i.e.}, 
$$ \lim_{i \rightarrow \infty} T_{li} \geq  \frac{1}{d z_{in}^*} \left[ \frac{d}{c} - \varepsilon \right]  > 0.  $$

\subsection{Proof of Lemma \ref{lem4}}

Let  $x_{o} \in \mathcal{D}_i$ and let $\hat{z}_{o2} = \frac{d/c}{2^{i-1}}$.  There is no loss of generality since if $x_{o} \not\in \mathcal{D}_i$ Assumption \ref{assd} trivially  holds over the $i$th cycle and if $\hat{z}_{o2} \neq \frac{d/c}{2^{i-1}}$, $\hat{z}_2(t,j) = \frac{d/c}{2^{i-1}}$ for some $t+j < \infty$.  Moreover,  the following reasoning applies {\it mutatis mutandis} if $\hat{z}_{o2} = -\frac{d/c}{2^{i-1}}$. We also use the fact that $|\tilde{z}_{o}| \leq g(i-1) \frac{\varepsilon}{\gamma}$ together with  \eqref{eqSobs} to conclude that
\begin{align} \label{eqztild}
|\tilde{z}(t,j)| \leq \varepsilon g(i-1) \qquad \forall (t,j) \in \mathop{\rm dom}\nolimits \tilde{z}.  
\end{align}

\begin{enumerate}[wide = 0pt, leftmargin = 1.3em,itemsep=3pt]
\item  At this point,  we estimate a lower bound on the  flow time that $\hat{z}_2$ takes to flow from $\frac{d/c}{2^{i-1}}$  to  $\frac{-d/c}{2^{i-1}}$.   Using \eqref{eqztild},  we conclude that such a time is lower bounded by the time $z_2$ takes to flow from $\frac{d/c}{2^{i-1}} - \frac{\varepsilon}{2^{i-1}}$ to $\frac{-d/c}{2^{i-1}} + \frac{\varepsilon}{2^{i-1}}$ when 
$$z_1 = - \left( \frac{z_{in}^*}{2^{i-1}} +  \frac{z_{in}^*}{2^{i-1}} g(i-1) \right).$$
 Let us denote such time by $T_{li}$,  which can be easily obtained by solving the ordinary differential equation
$$  \dot{z}_2 = - \left( \frac{ z_{in}^*}{2^{i-1}} +  \frac{z_{in}^*}{2^{i-1}} g(i-1) \right) (cz_2 + d), $$
with $z_{2o} = \frac{d/c - \varepsilon}{2^{i-1}}$.  After Item (ii) of  Lemma \ref{lem2-3} there exists $T_{lmin} > 0$ such that $T_{li} \geq T_{lmin}$ for all $i \in \{1,2,... \}$.  

\item  During the phase when  $\hat{z}_2$ flows from $\frac{d/c}{2^{i-1}}$  to  $\frac{-d/c}{2^{i-1}}$,  
  $z^* =  -\frac{z_{in}^*}{2^i}$.  Next, we show how to choose $k'>0$ to conclude that $z_{1e}$ must take at most $\frac{T_{li}}{2}$ units of time to enter the ball of radius $\frac{|z^*|}{2}$.  To  that end,  we use the Lyapunov function 
$$v (z_{1e}) := z_{1e}^2,  $$ 
whose time derivative along the solutions to $\mathcal{H}_i$ satisfies 
\begin{align*}
  \dot{v}  \leq  - 2 k'  v +   \frac{2 a  z_{in}^*}{2^{i}} \varepsilon g(i-1)   
  |z_{1e}|.
\end{align*}
By assuming, for example that $k' \geq 1$,  we conclude that
\begin{align*}
\dot{v} \leq -  k'  v +   \frac{a^2 \varepsilon^2 z_{in}^{*2}}{2^{2i}} g(i-1)^2,
\end{align*}
 Then,  for each $t \geq 0$ such that  $(t,0) \in \mathop{\rm dom}\nolimits x$,  
\begin{align*} 
v(t,0)  \leq  v(0,0) e^{-k't}   + \left[ \frac{a \varepsilon z_{in}^*}{2^{i} \sqrt{k'}} g(i-1) \right]^2  \!\! \big[ 1 - e^{-k't}  \big]
\end{align*}
so, choosing $k'$  large such that $k' \geq 4 a^2 \varepsilon^2$, we obtain 
\begin{align*} 
 v(t,0) & \leq   v(0,0) e^{-k't}   + \left[ \frac{z_{in}^*}{2^{i+1}} g(i-1) \right]^2   \left[ 1 - e^{-k't}  \right]  \\ &
\leq \left[ \frac{2 z_{in}^*}{2^{i}} + \frac{2 z_{in}^*}{2^{i}} g(i-1) \right]^2  e^{-k't} \\ & + \left[ \frac{z_{in}^*}{2^{i+1}} g(i-1) \right]^2   \left[ 1 - e^{-k't}  \right]
\\ & \leq 9 z^{*2} e^{-k't}   
 + \frac{z^{*2}}{4} \left[ 1 - e^{-k't}  \right].
\end{align*}
Next,  by letting $k' \geq - \frac{2 \ln (2^{-5})} {T_{lmin}}$,  we conclude that
\begin{align} \label{eqz1eb}
z_{1e}(t,0)^2 \leq \left( \frac{z^{*}}{2}\right)^2 \qquad \forall t \in \left[ \frac{T_{li}}{2}, T_{li} \right].  
\end{align}
Hence, during the interval where $\hat{z}_2$ flows from $\frac{d/c}{2^{i-1}}$ to $-\frac{d/c}{2^{i-1}}$,  we have $z_1 \in [z^*, z^*/2]$ for 
all $t$ belonging to a sub interval of length larger than $T_{li}/2$.  

\item  Next,  we estimate an upper bound on the time that 
$\hat{z}_2$ takes to flow from $\frac{d/c}{2^{i-1}}$  to  $-\frac{d/c}{2^{i-1}}$.    
Using \eqref{eqztild} and \eqref{eqz1eb},  we conclude that such a time is upper bounded by the time that $z_2$ takes to flow from $\frac{d/c}{2^{i-1}} + \frac{\varepsilon}{2^{i-1}}$ to $-\frac{d/c}{2^{i-1}} - \frac{\varepsilon}{2^{i-1}}$ when 
$$ z_1(t) = 
 \left\{ 
\begin{array}{ll} 
 \frac{z_{in}^*}{2^{i-1}} + \frac{z_{in}^*}{2^{i-1}} g(i-1)  & \forall \, t \in \left[0,\frac{T_{li}}{2} \right]
\\[4pt]
- \frac{|z^*|}{2}   & \forall \, t \geq  \frac{T_{li}}{2}.
\end{array}
\right.   $$   
Let us denote such a time by $T_{ui}$. Note that $T_{ui}$ can be easily obtained by solving the linear switched dynamics:
\begin{align*} 
 \dot{z}_2   =
\left\{ \begin{array}{ll}
  \left[  \frac{z_{in}^*}{2^{i-1}} + \frac{z_{in}^*}{2^{i-1}} g(i-1) \right]  (cz_2 + d)  & \forall t \in [0, \frac{T_{li}}{2}]
\\[4pt]
 - \frac{ |z^*|}{2}  (cz_2 + d) &  \forall t \geq \frac{T_{li}}{2},
 \end{array}
 \right.
 \end{align*}
from the initial condition $z_{2o} = \frac{d/c + \varepsilon}{2^{i-1}}$.
\end{enumerate}

\vskip 40pt
Finally,  we conclude that Assumption \ref{assd} holds on Cycle $i$ with 
\begin{align*} 
\tau_{di} & := \frac{T_{li}}{2},  & \tau_{si} & := T_{ui} - \frac{T_{li}}{2}, 
\\
\bar{z}_i & :=\frac{z_{in}^*}{2^{i-1}} + \frac{z_{in}^*}{2^{i-1}} g(i-1),  &  \underline{z}_i & := \frac{|z^*|}{2},
\end{align*}
which completes the proof. \hfill\small $\blacksquare$

%---
%% \bibliographystyle{IEEEtran}
%% \bibliography{hybrid_abs}
% Generated by IEEEtran.bst, version: 1.14 (2015/08/26)

%--

\end{document}